\documentclass{amsart}

\newtheorem{theorem}{Theorem}[section]

\theoremstyle{definition}
\newtheorem{definition}[theorem]{Definition}

\theoremstyle{remark}

\usepackage{graphicx}
\usepackage{color}
\usepackage{amsmath}
\usepackage{amsfonts}
\usepackage{amssymb}
\newcommand{\be}{\begin{equation}}
\newcommand{\ee}{\end{equation}}

\newcommand{\ba}{\begin{array}}

\newcommand{\ea}{\end{array}}

\newcommand{\beq}{\begin{eqnarray}}

\newcommand{\eeq}{\end{eqnarray}}

\newtheorem{lm}{lemma}

\newtheorem{thee}{theorem}

\newtheorem{proo}{proposition}

\newtheorem{co}{corollary}

\newtheorem{rem}{remark}

\newtheorem{deff}{definition}

\newcommand{\bd}{\begin{deff}}

\newcommand{\ed}{\end{deff}}

\newcommand{\bl}{\begin{lm}}

\newcommand{\el}{\end{lm}}

\newcommand{\bp}{\begin{proo}}

\newcommand{\ep}{\end{proo}}

\newcommand{\bt}{\begin{thee}}

\newcommand{\et}{\end{thee}}

\newcommand{\bc}{\begin{co}}

\newcommand{\ec}{\end{co}}

\newcommand{\brm}{\begin{rem}}

\newcommand{\erm}{\end{rem}}

\usepackage{t1enc}

\newcommand{\newc}{\newcommand}

\let\ccdot\cdot
\def\cdot{\hbox to 2.5pt{\hss$\ccdot$\hss}}

\newc{\aR}{\mbox{\boldmath{$ R$}}}
\newc{\aS}{\mbox{\boldmath{$ S$}}}
\newc{\aT}{\mbox{\boldmath{$ T$}}}
\newc{\aW}{\mbox{\boldmath{$ W$}}}

\newc{\aK}{\mbox{\boldmath{$ K$}}}
\newc{\aL}{\mbox{\boldmath{$ L$}}}

\newcommand{\bbC}{\mathbb{C}}

\newcommand{\bbK}{\mathbb{K}}

\usepackage{amssymb}
\usepackage{amscd}

\newc{\obstrn}[2]{B^{#1}_{#2}}

\newcommand{\rpl}                         
{\mbox{$
\begin{picture}(12.7,8)(-.5,-1)
\put(0,0.2){$+$}
\put(4.2,2.8){\oval(8,8)[r]}
\end{picture}$}}

\newcommand{\lpl}                         
{\mbox{$
\begin{picture}(12.7,8)(-.5,-1)
\put(2,0.2){$+$}
\put(6.2,2.8){\oval(8,8)[l]}
\end{picture}$}}

\usepackage{ifthen}

\newcommand{\bbR}{\mathbb{R}}

\newcommand{\og}{\mathbf{O}}

\newc{\tensor}[1]{#1}
\newc{\Mvariable}[1]{\mbox{#1}}
\newc{\down}[1]{{}_{#1}}
\newc{\up}[1]{{}^{#1}}

\newc{\JulyStrut}{\rule{0mm}{6mm}}
\newc{\midtenPan}{\mbox{\sf S}}
\newc{\midten}{\mbox{\sf T}}
\newc{\midtenEi}{\mbox{\sf U}}
\newc{\ATen}{\mbox{\sf E}}
\newc{\BTen}{\mbox{\sf F}}
\newc{\CTen}{\mbox{\sf G}}

\def\sideremark#1{\ifvmode\leavevmode\fi\vadjust{\vbox to0pt{\vss
 \hbox to 0pt{\hskip\hsize\hskip1em
 \vbox{\hsize3cm\tiny\raggedright\pretolerance10000
 \noindent #1\hfill}\hss}\vbox to8pt{\vfil}\vss}}}%

\numberwithin{equation}{section}

\newcounter{romenumi}
\newcommand{\labelromenumi}{(\roman{romenumi})}

\newcommand{\ten}{\Upsilon}

\begin{document}
\title{Deforming a Lie algebra by means of a two form}

\author{Pawe\l~ Nurowski} \address{Instytut Fizyki Teoretycznej,
Uniwersytet Warszawski, ul. Hoza 69, Warszawa, Poland}
\email{nurowski@fuw.edu.pl} \thanks{This research was supported by
the KBN grant 1 P03B 07529}

\date{\today}

\begin{abstract} We consider a vector space $V$ over $\bbK=\bbR$ or $\bbC$, equipped with a skew 
symmetric bracket $[\cdot,\cdot]: V\times V\to V$ and a 2-form 
$\omega:V\times V\to \bbK$.  A simple change of the Jacobi identity to 
the form $[A,[B,C]]+[C,[A,B]]+[B,[C,A]]=\omega(B,C)A+\omega(A,B)C+\omega(C,A)B$ opens 
new possibilities, which shed new light on the Bianchi classification of 
3-dimensional Lie algebras.

\end{abstract}
\maketitle
\section{Introduction}
In reference \cite{bobi} we considered a real vector space $V$ of dimension $n$ equipped with a Riemannian metric $g$ and a symmetric 3-tensor $\ten_{ijk}$ such that: i) $\ten_{ijk}=\ten_{(ijk)}$, ii) $\ten_{ijj}=0$ and iii) $\ten_{jki}\ten_{lmi}+\ten_{lji}\ten_{kmi}+\ten_{kli}\ten_{jmi}=g_{jk}g_{lm}+g_{lj}g_{km}+g_{kl}g_{jm}.$ Such tensor defines a bilinear product 
$\{\cdot,\cdot\}:V\times V\to V$ given by $$\{A,B\}_i=\ten_{ijk}A_jB_k.$$ 
This product is symmetric 
\be\{A,B\}=\{B,A\}\label{sym}\ee 
due to property ii), 
and it satisfies a three-linear identity:
\be
\{A,\{B,C\}\}+\{C,\{A,B\}\}+\{B,\{C,A\}\}=g(B,C)A+g(A,B)C+g(C,A)B,
\label{jac}\ee
due to property iii). Restricting our attention to 
structures $(V,g,\{\cdot,\cdot\})$ associated with tensors $\ten$ as above, we note that they are related to the isoparametric hypersurfaces in spheres  \cite{bryant,Cartan1}. Using Cartan's results \cite{cartan} on isoparametric hypersurfaces we concluded in \cite{pn} that structures $(V,g,\{\cdot,\cdot\})$ exist only in dimensions $5,8,14$ and $26$.

A striking feature of property (\ref{jac}) is that it resembles very much the Jacobi identity satisfied by every Lie algebra. The main difference is that for a Lie algebra the bracket $\{\cdot,\cdot\}$ should be {\it anti}-symmetric and that the analog of (\ref{jac}) should have r.h.s equal to {\it zero}. 

Adapting properties (\ref{sym})-(\ref{jac}) to the notion of a Lie algebra we are led to the following structure.
\begin{definition}
A vector space $V$ equipped with a bilinear bracket $[\cdot,\cdot ]:V\times V\to V$ and a 2-form $\omega: V\times V\to \bbK=\bbR$ or $\bbC$ such that 
$$[A,B]=-[B,A]\quad\quad\quad {\rm and}$$
\be [A,[B,C]]+[C,[A,B]]+[B,[C,A]]=\omega(B,C)A+\omega(A,B)C+\omega(C,A)B
\label{jacc}
\ee
is called an $\omega$-deformed Lie algebra.  
\end{definition} 
This definition obviously generalizes the notion of a Lie algebra and coincides with it when $\omega\equiv 0$. Note also that if the dimension of $V$ is $\dim V=2$, then $\omega(B,C)A+\omega(A,B)C+\omega(C,A)B\equiv 0$ for any 2-form $\omega$ and every $A,B,C\in V$. Thus in 2-dimensions it is impossible to $\omega$-deform the Jacobi identity, and 2-dimensional $\omega$-deformed Lie algebras are just the Lie algebras equipped with a 2-form $\omega$. This is not anymore true 
if $\dim V\neq 2$. Indeed assuming that $\dim V\neq 2$ and that $\omega(B,C)A+\omega(A,B)C+\omega(C,A)B\equiv 0$ for all $A,B,C\in V$ we easily prove that $\omega\equiv 0$.  

The aim of this note is 
to show that there exist $\omega$-deformed Lie algebras in dimensions greater than 2 which are not just the Lie algebras. 
\section{Dimension 3.}
It follows that if $\dim V\leq 2$ then all the $\omega$-deformed Lie algebras are just the Lie algebras. To show that in $\dim V=3$ the situation is different we follow the procedure used in the Bianchi classification \cite{bian} of 3-dimensional Lie algebras.

Let $\{e_i\}$, $i=1,2,3,$ be a basis of an $\omega$-deformed 3-dimensional Lie algebra. Then, due to skew-symmetry, we have $[e_i,e_j]=c^k_{~ij}e_k$, $\omega(e_i,e_j)=\omega_{ij}$, where $c^k_{~ij}=-c^k_{ji}$ and 
$\omega_{ij}=-\omega_{ji}$. Due to the $\omega$-deformed Jacobi identity (\ref{jacc}), we also have  
$$
c^m_{~li}c^i_{~jk}+c^m_{~ki}c^i_{~lj}+c^m_{~ji}c^i_{~kl}=\delta^m_{~l}\omega_{jk}+\delta^m_{~k}\omega_{lj}+\delta^m_{~j}\omega_{kl},
$$
which is equivalent to
\be
c^m_{~i[l}c^i_{~jk]}+\delta^m_{~[l}\omega_{jk]}=0.
\label{jaC}\ee

We now find all the orbits of the above defined pair of tensors $(c^k_{~ij},\omega_{ij})$ under the action of the group ${\bf GL}(3,\bbR)$.

We recall that in three dimensions, we have the totally skew symmetric 
Levi-Civita symbol $\epsilon_{ijk}$, and its totally skew symmetric inverse 
$\epsilon^{ijk}$ such that $\epsilon_{ijk}\epsilon^{ilm}=\delta^l_{~j}\delta^m_{~k}-\delta^m_{~j}\delta^l_{~k}$.  This can be used to rewrite the $\omega$-deformed Jacobi identity (\ref{jaC}). Indeed, since in 3-dimensions every totally skew symmetric 3-tensor is proportional to $\epsilon_{ijk}$, the l.h.s. of (\ref{jaC}) can be written as
$$t^m=0\quad\quad\quad{\rm with}\quad\quad\quad t^m=(c^m_{~il}c^i_{~jk}+\delta^m_{~l}\omega_{jk})\epsilon^{ljk}.$$
In addidtion, we may use $\epsilon_{ijk}$ to write $c^i_{~jk}$ as
\be
c^i_{~jk}=n^{il}\epsilon_{jkl}-\delta^i_{~j}a_k+\delta^i_{~k}a_j,\label{pier}
\ee  
where the symmetric matrix $n^{il}$ is related to $c^k_{~ij}$ via
$$ n^{il}=\tfrac12 (c^{il}+c^{li}),\quad\quad\quad{\rm with}\quad\quad\quad c^{il}=\tfrac12 c^i_{~jk}\epsilon^{jkl}.$$ The vector $a_m$ is related to $c^k_{~ij}$ via $$a_m=\tfrac12 \epsilon_{mil}c^{il}.$$ 
Similarly, we write $\omega_{ij}$ as
\be
\omega_{ij}=\epsilon_{ijk}b^k,\label{dru}
\ee
with 
$$b^k=\tfrac12 \epsilon^{mik}\omega_{ik}.$$
Thus, in three dimensions the structural constants $(c^k_{~ij},\omega_{ij})$ of the $\omega$-deformed Lie algebra are uniquely determined via (\ref{pier}), (\ref{dru}) by specifying a symmetric matrix $n^{il}$ and two vectors $a_m$ and $b^k$. In terms of the triple $(n^{il},a_m,b^k)$ the vector $t^m$ is given by $t^m=4n^{ml}a_l+2 b^m$, so that the $\omega$-deformed Jacobi identity 
(\ref{jaC}) is simply
\be
b^i=-2n^{il}a_l. \label{JaC}
\ee  
Thus, given $n^{il}$ and $a_m$, the vector $b^m$ defining $\omega$ is totally determined.
Now we use the action of ${\bf GL}(3,\bbR)$ group to bring $n^{il}$ to the diagonal form (it is always possible since $n^{il}$ is symmetric), so that 
$$n^{il}={\rm diag}(n^1,n^2,n^3).$$
It is obvious that without loss of generality we always can have
$$n^i=\pm 1, 0\quad\quad\quad i=1,2,3.$$
After achiving this we may still use an orthogonal transformation preserving the matrix $n^{il}$ to bring the vector $a_m$ to a simpler form then $a_m=(a_1,a_2,a_3)$. For example in the case $n^{il}={\rm diag}(1,1,1)$ we may always achieve $a_m=(0,0,a)$.
Thus to represent a ${\bf GL}(3,\bbR)$ orbit of $(c^i_{~jk},\omega_{ij})$ it is enough to take $n^{il}$ in the diagonal form with the diagonal elements being equal to $\pm 1,0$ and to take $a_m$ in the simplest possible form obtainable by the action of $\og(n^{il})$. Finally we notice that the 
so specified choice of $n^{il}$ is still preserved when the basis is scalled according to
\be
e_1\to \lambda_1 e_1,\quad\quad e_2\to \lambda_2 e_2\quad\quad e_3\to \lambda_3 e_3,\label{resfr}\ee 
with  
$$(\lambda_1\lambda_2-\lambda_3)n_3=0,\quad\quad (\lambda_3\lambda_1-\lambda_2)n_2=0,\quad\quad(\lambda_2\lambda_3-\lambda_1)n_1=0.$$
These transformations can be used to scale the vector $a_m$ via 
$$a_m\to (\lambda_1 a_1,\lambda_2 a_2,\lambda_3 a_3).$$ 

We are now in a position to give the full classification of 3-dimensional $\omega$-deformed Lie algebras. In all the types of the classification the commutation relations and the $\omega$ are given by:
\begin{eqnarray*}&&[e_1,e_2]=n^3 e_3-a_2 e_1+a_1 e_2,\quad\quad [e_3,e_1]=n^2 e_2-a_1 e_3+a_3 e_1,\\
&&[e_2,e_3]=n^1 e_1-a_3 e_2+a_2 e_3
\end{eqnarray*}
$$\omega(e_1,e_2)=-2n^3 a_3,\quad\quad \omega(e_3,e_1)=-2n^2a_2,\quad\quad \omega(e_2,e_3)=-2n^1a_1.$$
The classification splits into two main branches depending on vanishing or not of $a_m$.

If $a_m=0$, then $b^m=0$ and all the possibilities are given in the following table:\begin{center}
{\small 
\begin{tabular}{|l|l|l|l|}
\hline
Bianchi type& $n^1$&$n^2$&$n^3$\\
\hline\hline 
$I$&0&0&0\\
$II$&1&0&0\\
$VI_0$&1&-1&0\\
$VII_0$&1&1&0\\
$VIII$&1&1&-1\\
$IX$&1&1&1\\
\hline
\end{tabular}\hspace{1cm} $a_m=0$, $b^m=0$.
}
\end{center}
All types from this table have $\omega=0$ and as such 
correspond to the usual 3-dimensional Lie algebras.

If $a_m\neq 0$ then, depending on the signature of $n^{il}$, vector $a_m$ may be spacelike, timelike, null or degenerate. The orthogonal transforamtions we use to normalize this vector preserve its type, so the classification splits according to the causal properties of $a_m$. If $n^2=n^3=0$ or $n^1=-n^2=1$, $n^3=0$, we may use transformations (\ref{resfr}) to totally fix $a_m$. This leads to types $V$, $IV$, $IV_x$ and $VI_x, VI_y, VI_n$ below. In all other cases transformations (\ref{resfr}) can be used to express $a_m$ in terms of only one parameter $a>0$ so that the different positive parameters $a$ correspond to nonequivalent algebras. The resulting classification is summarized in the following table:
\begin{center}
{\small
\begin{tabular}{|l|l|l|l|l|l|}
\hline
Bianchi type& $n^1$&$n^2$&$n^3$&$a_m$&$b^m$\\
\hline\hline 
$V$&0&0&0&(0,0,1)&(0,0,0)\\
\hline
$IV$&1&0&0&(0,0,1)&(0,0,0)\\
$IV_x$&1&0&0&(1,0,0)&(-2,0,0)\\
\hline
$VI_a$&1&-1&0&$(0,0,a>0)$&(0,0,0)\\
$VI_{x}$&1&-1&0&$(1,0,0)$&(-2,0,0)\\
$VI_{y}$&1&-1&0&$(0,1,0)$&(0,2,0)\\
$VI_{n}$&1&-1&0&$(1,1,0)$&(-2,2,0)\\
\hline
$VII_a$&1&1&0&$(0,0,a>0)$&(0,0,0)\\
$VII_{x}$&1&1&0&$(1,0,0)$&(-2,0,0)\\
\hline
$VIII_a$&1&1&-1&$(0,0,a>0)$&$(0,0,2a)$\\
$VIII_{xa}$&1&1&-1&$(a>0,0,0)$&$(-2a,0,0)$\\
$VIII_{na}$&1&1&-1&$(a>0,0,a)$&$(-2a,0,2a)$\\
\hline
$IX_a$&1&1&1&$(0,0,a>0)$&$(0,0,-2a)$\\
\hline
\end{tabular}.
}
\end{center}
In the above two tables all the types which have $b^m=0$ are just the usual 3-dimensional Lie algebras. Apart from the types I and V all the Bianchi types admit $\omega$ deformation. It is interesting to note that types $VIII$ and $IX$, which in the Lie algebra setting do not admit $a_m\neq 0$ deformation, admit a one-parameter $\omega$-deformations. 

We have the following theorem.
\begin{theorem}
All the 3-dimensional $\omega$-deformed Lie algebras are given in the following table
\begin{center}
{\small
\begin{tabular}{|l|l|l|l|l|l|}
\hline
Bianchi type& $n^1$&$n^2$&$n^3$&$(a_1,a_2,a_3)$&$(b^1,b^2,b^3)$\\
\hline\hline 
$IV_x$&1&0&0&(1,0,0)&(-2,0,0)\\
\hline
$VI_{x}$&1&-1&0&$(1,0,0)$&(-2,0,0)\\
$VI_{y}$&1&-1&0&$(0,1,0)$&(0,2,0)\\
$VI_{n}$&1&-1&0&$(1,1,0)$&(-2,2,0)\\
\hline
$VII_{x}$&1&1&0&$(1,0,0)$&(-2,0,0)\\
\hline
$VIII_a$&1&1&-1&$(0,0,a>0)$&$(0,0,2a)$\\
$VIII_{xa}$&1&1&-1&$(a>0,0,0)$&$(-2a,0,0)$\\
$VIII_{na}$&1&1&-1&$(a>0,0,a)$&$(-2a,0,2a)$\\
\hline
$IX_a$&1&1&1&$(0,0,a>0)$&$(0,0,-2a)$\\
\hline
\end{tabular}
}
\end{center}
They satisfy the commutation relations 
\begin{eqnarray*}&&[e_1,e_2]=n^3 e_3-a_2 e_1+a_1 e_2,\quad\quad [e_3,e_1]=n^2 e_2-a_1 e_3+a_3 e_1,\\
&&[e_2,e_3]=n^1 e_1-a_3 e_2+a_2 e_3
\end{eqnarray*}
$$\omega(e_1,e_2)=-2n^3 a_3,\quad\quad \omega(e_3,e_1)=-2n^2a_2,\quad\quad \omega(e_2,e_3)=-2n^1a_1.$$
with the real parameters $(n^1,n^2,n^3,a_1,a_2,a_3)$ specified in the table. Algebras corresponding to different $(n^1,n^2,n^3,a_1,a_2,a_3)$ are nonequivalent.
\end{theorem}  

Finally we show that any $\omega$-deformed Lie algebra must have quite nontrivial structure constants. Indeed, in any dimension $\dim V=n>2$ the structure constants of an $\omega$-deformed Lie algebra, which are defined by $[e_i,e_j]=c^k_{~ij}e_k$, may be decomposed as follows:
$$c^i_{~jk}=\alpha^i_{~jk}+a_k\delta^i_{~j}-a_j\delta^i_{~k},$$
where 
$$\alpha^i_{~ik}=0,\quad\quad\quad a_k=\tfrac{1}{n-1}c^i_{~ik}.$$
Then a simple calculation using the $\omega$-deformed Jacobi identity (\ref{jacc}) shows that
$$\omega(e_j,e_k)= \tfrac{n-1}{n-2}a_i\alpha^i_{~jk}.$$
This shows that nonvanishing $\omega$ is only possible if both $a_i$ and $\alpha^i_{~jk}$ are nonvanishing. 
\section{Acknowledgements}
I am very grateful to Jose Figueroa-O'Farrill for reading the draft of this paper and correcting an error in my enumeration of the Bianchi types. I also wish 
to thank David Calderbank for helpful discusions.

\end{document}